\documentclass[12pt,scheme=plain,fontset=none]{ctexart}
\usepackage[T1]{fontenc}
\usepackage[utf8]{inputenc}
\usepackage{lmodern}

\usepackage[a4paper,margin=1in]{geometry}
\usepackage{amsmath,amssymb,amsthm,mathtools}
\usepackage{enumitem}
\usepackage{needspace}
\usepackage{microtype}
\usepackage{xcolor}
\usepackage{hyperref}
\usepackage{cite}

\hypersetup{
  colorlinks=true,
  linkcolor=blue,
  citecolor=blue,
  urlcolor=blue,
  pdftitle={Frame Measure Rigidity and Fiber Entropy for Homogeneous Self-Similar Measures},
  pdfauthor={Wei-Jie Wang},
  pdfsubject={Fourier frame rigidity, conditional coding measures, projection entropy, and Hausdorff dimension},
  pdfkeywords={Fourier frame rigidity, Fourier frame measure, homogeneous self-similar measure, conditional coding measure, projection entropy, dimension rigidity}
}

\numberwithin{equation}{section}
\AtBeginDocument{%
  \setlength{\abovedisplayskip}{6.5pt plus 2pt minus 2pt}%
  \setlength{\belowdisplayskip}{6.5pt plus 2pt minus 2pt}%
  \setlength{\abovedisplayshortskip}{4pt plus 2pt minus 2pt}%
  \setlength{\belowdisplayshortskip}{4pt plus 2pt minus 2pt}%
}

\theoremstyle{plain}
\newtheorem{theorem}{Theorem}[section]
\newtheorem{corollary}[theorem]{Corollary}
\newtheorem{proposition}[theorem]{Proposition}
\newtheorem{lemma}[theorem]{Lemma}
\theoremstyle{definition}

\newtheorem{remark}[theorem]{Remark}

\newcommand{\R}{\mathbb R}
\newcommand{\N}{\mathbb N}
\newcommand{\Z}{\mathbb Z}
\newcommand{\Pp}{\mathbb P}
\newcommand{\E}{\mathbb E}
\newcommand{\D}{\mathcal D}
\newcommand{\B}{\mathcal B}
\newcommand{\G}{\mathcal G}
\newcommand{\cP}{\mathcal P}
\newcommand{\supp}{\operatorname{supp}}
\newcommand{\At}{\operatorname{At}}
\newcommand{\1}{\mathbf 1}
\newcommand{\dd}{\,\mathrm d}
\newcommand{\wh}{\widehat}
\newcommand{\norm}[1]{\lVert #1\rVert}
\newcommand{\abs}[1]{\lvert #1\rvert}
\newcommand{\dimH}{\operatorname{dim}_{\mathrm H}}

\title{Frame Measure Rigidity and Fiber Entropy for Homogeneous Self-Similar Measures}
\author{ Wei-Jie Wang}
\date{}

\begin{document}
\maketitle

\begin{abstract}
Let $N\ge2$, $0<\rho<1$, and let
$b_0,\ldots,b_{N-1}\in\R^d$ be distinct translation vectors. Let $\mu$
be the uniform homogeneous self-similar probability measure on $\R^d$
generated by $\{S_a(x)=\rho x+b_a\}_{a=0}^{N-1}$, with no separation
assumption. We prove that if $\mu$ admits a frame measure, then its
conditional coding measure is a Dirac mass for $\mu$-almost every point.
Equivalently, the coding map restricts to a Borel bijection between
full-measure subsets of the symbolic and physical spaces. It follows that
the mean fiber entropy vanishes at every finite level and that $\mu$ is
exact dimensional, with
\[
  \dimH\mu=\frac{\log N}{-\log\rho}\le d.
\]
Thus, if $\dimH\mu<\frac{\log N}{-\log\rho}$, then $\mu$ admits no
frame measure. In particular, no frame measure exists when $N\rho^d>1$.
As an application, the uniform consecutive-digit Bernoulli convolution
on $\R$ admits a frame measure for $N^{-1}\le\rho<1$ if and only if
$\rho=N^{-1}$.
\end{abstract}



\section{Introduction}\label{sec:intro}

Let $\mu$ be a compactly supported Borel probability measure on $\R^d$. For $\lambda\in\R^d$, we
write
\[
  e_\lambda(x):=e^{2\pi i\lambda\cdot x},\qquad x\in\R^d.
\]
A countable family $E(\Lambda):=\{e_\lambda:\lambda\in\Lambda\}$ is called a
\emph{Fourier frame} for $L^2(\mu)$ if there exist constants
$0<A\le B<\infty$ satisfying the frame inequality
\begin{equation}\label{eq:ordinary-frame}
  A\norm{f}_{L^2(\mu)}^2
  \le \sum_{\lambda\in\Lambda}
       \abs{\langle f,e_\lambda\rangle_{L^2(\mu)}}^2
  \le B\norm{f}_{L^2(\mu)}^2,
  \qquad \forall f\in L^2(\mu).
\end{equation}
If $E(\Lambda)$ forms an orthonormal basis for $L^2(\mu)$, then $\Lambda$ is called a
\emph{spectrum} of $\mu$ and $\mu$ is called a \emph{spectral measure}.
Thus, every spectrum of \(\mu\) gives rise to a Fourier frame for \(\mu\).
Spectrality has been studied for several classes of Cantor-type measures;
see, for example, the work of Dutkay, Han and Sun, Dutkay and Jorgensen,
and Laba and Wang~\cite{DutkayHanSun2009,DutkayJorgensen2007,LabaWang2002}.
For Bernoulli convolutions and consecutive-digit measures, related results
were obtained by Dai, Dai--He--Lai, and Dai--He--Lau
\cite{Dai2012,DaiHeLai2013,DaiHeLau2014}.

\vspace{0.2cm}

A frame measure allows the frequency measure in the Fourier-frame
inequality to be any positive Borel measure, rather than only a discrete
measure. For $f\in L^2(\mu)$, define the Fourier transform of the measure
$f\,d\mu$ by
\[
  \wh{f\mu}(\xi)
  :=\int_{\R^d}f(x)e^{-2\pi i\xi\cdot x}\dd\mu(x).
\]
A positive Borel measure $\Gamma$ on $\R^d$ is called a \emph{frame measure} for
$\mu$ if there exist constants
$0<A\le B<\infty$ such that
\begin{equation}\label{eq:frame-inequality}
  A\norm{f}_{L^2(\mu)}^2
  \le \int_{\R^d}\abs{\wh{f\mu}(\xi)}^2\dd\Gamma(\xi)
  \le B\norm{f}_{L^2(\mu)}^2,
  \qquad \forall f\in L^2(\mu).
\end{equation}
Taking $\Gamma=\sum_{\lambda\in\Lambda}\delta_\lambda$ recovers the classical discrete Fourier frame condition
\eqref{eq:ordinary-frame}. More generally, if $\Gamma=\sum_{\lambda\in\Lambda}c_{\lambda}\delta_\lambda$ with \(c_{\lambda}>0\), then \eqref{eq:frame-inequality} is precisely the frame inequality for the weighted Fourier frame $\{\sqrt{c_\lambda}e_\lambda:\lambda\in\Lambda\}$ associated with \(\mu\).
Thus, every discrete Fourier frame gives rise to an associated frame measure;  in turn, a bounded continuous frame may be sampled to produce a discrete one ~\cite{FreemanSpeegle2019}. In the exponential‑frame context, this sampling theorem permits non‑trivial frequency multiplicities.

\vspace{0.2cm}

We first recall some results on frame measures for self-similar measures.
Dutkay, Han and Weber~\cite{DutkayHanWeber2014} used frame measures to
study Fourier frames for singular self-similar measures, mainly under the
open set condition. Dutkay and Lai~\cite{DutkayLai2014} obtained a
necessary \textit{uniformity} condition for measures admitting frame
measures. For non-overlapping self-affine measures, this condition rules
out Fourier frames when the branch weights are unequal. Further constructions
of Fourier bases and weighted Fourier frames for self-affine and fractal
measures were given by Dutkay, Lai and Wang, and by Dutkay and Ranasinghe
\cite{DutkayLaiWang2017,DutkayRanasinghe2016}. Recent work has
also established nonexistence results for one-dimensional $N$-Bernoulli
convolutions under separation or arithmetic assumptions
\cite{FuSongWang2026,FuWang2026,PontLiehrTaylor2026}. Wang and
Tang~\cite{WangTang2026} proved that, for uniform self-similar measures
on $\R$ generated by any $N$ distinct real digits, positive conditional
entropy of the first digit given the coded point rules out a frame
measure. In particular, no such measure admits a frame measure when
$\rho>1/N$. Further constructions and
nonexistence results for Fourier frames of singular measures can be found
in~\cite{LaiWang2017,FuLai2018,Lev2018}.
 
%

 \vspace{0.2cm}
 
Most existing results on Fourier frames for self-similar measures use
separation assumptions or additional arithmetic conditions on the
underlying iterated function system (IFS). The problem is less understood
when the system has overlaps. In that case, the coding map from the
symbolic space to the attractor need not be injective, and the Hausdorff
dimension of the measure may be smaller than the similarity dimension.
Dimension results for overlapping self-similar systems can be found in
\cite{Hochman2014,Shmerkin2019,SimonToth2005,Solomyak1995}. Here we ask
whether the existence of a frame measure forces the coding map to be
almost everywhere injective, without a separation assumption. 
 


\vspace{0.2cm}

 For a topological space
$X$,  $\B(X)$ denotes its Borel $\sigma$-algebra. 
Fix an integer \(N\ge2\), a contraction ratio \(0<\rho<1\), and distinct translation vectors
\(
  b_0,\ldots,b_{N-1}\in\R^d.
\)
Define the associated homothetic iterated function system (IFS) by
\begin{equation}\label{eq:intro-ifs}
  S_a(x):=\rho x+b_a,
  \qquad a\in\D:=\{0,1,\ldots,N-1\}.
\end{equation}
This system determines a unique nonempty compact attractor
$K\subset\R^d$ satisfying
\begin{equation}\label{eq:intro-attractor}
	K=\bigcup_{a\in\D}S_a(K),
\end{equation}
and a unique probability measure \(\mu\) satisfying
\begin{equation}\label{eq:intro-invariance}
	\mu=\frac1N\sum_{a\in\D}(S_a)_\#\mu,
\end{equation}
where  $(S_a)_\#\mu$  denotes the pushforward measure of \(\mu\) under \(S_a\), given by
\begin{equation*}
	((S_a)_\#\mu)(E)=\mu(S_a^{-1}E), \quad \forall\ E\in\B(\R^d).
\end{equation*}
 Existence and uniqueness of \(K\) and \(\mu\) follow from the classic IFS fixed-point theorem of Hutchinson~\cite{Hutchinson1981}.  Then $K$ and \(\mu\) are referred to as the self‑similar set and self‑similar measure, respectively.

\vspace{0.2cm}

We give necessary conditions for the measure $\mu$ in
\eqref{eq:intro-invariance} to admit a frame measure. Our main result
shows that frame-measure existence forces almost-everywhere uniqueness
of coding and, consequently, rules out a dimension drop. We now introduce
the symbolic notation needed to state the result. 

\vspace{0.2cm}

Let \(\Omega:=\D^{\N}\) denote the one-sided infinite symbol space endowed with the product Bernoulli measure
\begin{eqnarray}\label{eq:product-Bernoulli measure}
  \Pp:=\left(\frac1N\sum_{a\in\D}\delta_a\right)^{\otimes\N}.
\end{eqnarray}
Define the coding map \(\pi: \Omega\to K\) by
\begin{equation*}
  \pi(w)
  :=\sum_{j=1}^{\infty}\rho^{j-1}b_{w_j}
  \in\R^d.
\end{equation*}
Then \(\mu=\pi_\#\Pp\). For each spatial point \(x\in K\), let $\{\Pp_x:x\in K\}$ denote the family of conditional probability distributions on \(\Omega\) supported on the symbolic addresses \(w\) satisfying \(\pi(w)=x\). 
For a finite word
$u=(u_1,\ldots,u_n)\in\D^n$, define its cylinder set by
\[
  [u]:=\{w\in\Omega:w_j=u_j\text{ for }1\le j\le n\}.
\]
Define the level-$n$ fiber entropy by
\begin{equation}\label{eq:intro-Rn}
  R_n
  :=\int_K\left[-\sum_{u\in\D^n}
       \Pp_x([u])\log\Pp_x([u])\right]\dd\mu(x).
\end{equation}
with the convention \( 0\log0:=0.\)
The equality \(R_n=0\) holds if and only if, for \(\mu\)-a.e. $x\in K$, the conditional coding distribution is concentrated on a single level-$n$ cylinder.

\vspace{0.2cm}

We adopt the standard definition of Hausdorff dimension for measures: 
\begin{equation}\label{eq:measure-dimension}
	\dimH\mu
	:=\inf\{\dimH E:E\in\B(\R^d),\ \mu(E)=1\}.
\end{equation}
The similarity dimension of the IFS \(\{S_a\}_{a\in\D}\) is \(\frac{\log N}{-\log\rho}\). This quantity may exceed $d$, whereas every Borel probability measure on $\R^d$ satisfies $\dimH\mu\le d$.

\begin{theorem}\label{thm:rigidity}
Let $b_0,\ldots,b_{N-1}\in\R^d$ be distinct translation
vectors, and let $\mu$ be the uniform homogeneous self-similar measure
defined by \eqref{eq:intro-invariance}.  If $\mu$ admits a frame measure on
$\R^d$, then the following statements hold.
\begin{enumerate}[label=\textup{(\roman*)},leftmargin=2.2em]
\item There exist Borel sets $\Omega_0\subset\Omega$ and $K_0\subset K$ of
full $\Pp$- and $\mu$-measure, respectively, such that the coding restriction
\(
  \pi:\Omega_0\longrightarrow K_0
\)
is a Borel bijection.  Equivalently, the conditional coding measure $\Pp_x$ is a Dirac mass for
$\mu$-almost every $x\in K$.
\item The  level-\(n\) fiber entropy satisfies
\(
  R_n=0
\)
for every integer $n\ge1$.
\item The measure $\mu$ is exact dimensional, with Hausdorff dimension
\(
  \dimH\mu=\frac{\log N}{-\log\rho}\le d.
\)
\end{enumerate}
\end{theorem}

The conditions in Theorem~\ref{thm:rigidity} are necessary but not sufficient for the existence of a frame measure. For example, the middle-third Cantor measure has essentially unique coding but admits no frame measure~\cite{FuWang2026}.

\begin{remark}
Lai~\cite{Lai2011} completely characterized compactly supported absolutely continuous measures admitting Fourier frames: a Fourier frame exists for such a measure if and only if its density is essentially bounded from above and below on its support. In particular, for equal-weight self-similar measures $\mu$ defined in \eqref{eq:intro-invariance} on \(\mathbb R\), Theorem~\ref{thm:rigidity}\textup{(iii)} shows that absolute continuity together with Fourier frame existence forces the critical contraction ratio \(\rho=1/N\). On the other hand, $\mu_{\rho,N}$ is known to be absolutely continuous for almost every \(\rho\in(1/N,1)\) (see~\cite{ShmerkinSolomyak2016, SimonToth2005, Solomyak1995}). Combined with Lai's theorem~\cite{Lai2011}, this rules out Fourier frames of $\mu_{\rho,N}$ for all \(\rho\in(1/N,1)\) up to a null exceptional set. Theorem \ref{thm:rigidity} removes this null exceptional set entirely: for every \(\rho>1/N\), the measure $\mu_{\rho,N}$ admits no frame measure, regardless of absolute continuity or singularity.
\end{remark}

\begin{remark}
 Condition (i) in Theorem~\ref{thm:rigidity} implies that, for any two distinct words \(u,v\in\D^n\), the pushforward measures \((S_u)_\#\mu\) and \((S_v)_\#\mu\) are mutually singular; see \eqref{eq:exact-cylinder-measure}. This does not by itself assert that \(\mu(S_u(K)\cap S_v(K))=0\), nor does it imply that the geometric cylinders are disjoint or that the open set condition holds. Related measure-theoretic separation results for spectral self-affine measures were obtained by Deng and Chen~\cite{DengChen2021}. 
\end{remark}

\begin{remark}
	The dimension identity (iii) in Theorem~\ref{thm:rigidity} provides a convenient criterion for ruling out frame measures whenever the similarity dimension exceeds the ambient dimension, or equivalently, whenever $N\rho^d>1$.
\end{remark}

If \(\rho^d>1/N\), then \(-\frac{\log N}{\log\rho}>d\). Since every Borel probability measure on \(\mathbb{R}^d\) satisfies \(\dim_H\mu\le d\), Theorem~\ref{thm:rigidity} gives the following corollary.




\begin{corollary}\label{cor:supercritical}
	If $\rho^d>1/N$, then $\mu$ admits no frame measure.  
\end{corollary}

For $b_a=\rho a$, the measure in
\eqref{eq:intro-invariance} is the uniform $N$-Bernoulli
convolution
\begin{equation}\label{eq:intro-bernoulli}
	\mu_{\rho,N}
	:=\mathop{\ast}_{j=1}^{\infty}
	\left(\frac1N\sum_{a=0}^{N-1}\delta_{\rho^ja}\right).
\end{equation}
For this family,
Corollary~\ref{cor:supercritical} combined with the critical case
$\rho=N^{-1}$, gives a complete classification over the range
$N^{-1}\le\rho<1$: 
a frame measure exists for \(\mu_{\rho,N}\) \(\Longleftrightarrow\)
$\rho=N^{-1}$. 
At the critical value,
\begin{equation}\label{eq:intro-critical-measure}
	\mu_{1/N,N}=\1_{[0, 1]}\dd x,
\end{equation}
and $\{e^{2\pi ikx}:k\in\Z\}$ forms an orthonormal basis of
$L^2(\mu_{1/N,N})$. We return to this application in
Section~\ref{sec:further-questions}.

The main step is to prove part~\textup{(i)}. Once part~\textup{(i)} is established, the remaining two conclusions follow as
\[
  \Pp_x\text{ is a Dirac measure for }\mu\text{-a.e. }x
  \quad\Longrightarrow\quad
  R_n=0,\ \forall n\ge1
  \quad\overset{\eqref{eq:feng-hu-specialized}}{\Longrightarrow}\quad
  \dimH\mu=\frac{\log N}{-\log\rho}.
\]

\vspace{0.2cm}

When the coding map is not one-to-one, the conditional measures $\Pp_x$ on the fibers of $\pi$ need not be Dirac. The associated finite-level conditional entropies are related to projection entropy. Feng and Hu~\cite{FengHu2009} developed the projection-entropy framework and proved the exact-dimensionality theorem used below. Related dimension and overlap results for self-similar measures were obtained by Hochman~\cite{Hochman2014} and Shmerkin~\cite{Shmerkin2019}. In the Bernoulli-convolution setting, almost-everywhere absolute-continuity results go back to Solomyak~\cite{Solomyak1995} and, for consecutive-digit systems, Simon and T\'oth~\cite{SimonToth2005}. We prove that if $\mu$ admits a frame measure, then $\Pp_x$ is a Dirac mass for $\mu$-almost every $x$. Consequently, $R_n=0$ for every $n\ge1$, and the Feng--Hu dimension formula gives $\dimH\mu=\frac{\log N}{-\log\rho}$.

\vspace{0.2cm}

The paper is organized as follows. In Section~\ref{sec:posterior-frame}, we set up symbolic dynamical and measure-theoretic preliminaries, the posterior Fourier identity, and the equal-level comparison technique deployed throughout the subsequent analysis.
In
Section~\ref{sec:atomicity}, we derive the one-dimensional collision bound via the rightmost branch construction and establish pure atomicity of the conditional coding measures.
In Section~\ref{sec:entropy}, we prove the Dirac structure of conditional fibers and essential injectivity, and then apply the Feng--Hu dimension theorem using the first-level fiber entropy $R_1$.
In Section~\ref{sec:Rd}, we extend our main result to arbitrary dimensions. This is accomplished by constructing an explicit visible branch via an exposed translation and implementing the prior symbolic-entropic line of reasoning. In
Section~\ref{sec:further-questions}, we conclude the paper with discussions on practical applications, limitations of our necessary conditions, and related open problems for future investigation.

\section{Preliminaries}\label{sec:posterior-frame}

This section assembles all symbolic, measure-theoretic, and entropic notation required for subsequent proofs. We start with symbolic dynamics, disintegration theory, and  word maps associated with affine IFS, then introduce cylinder partitions, Shannon entropy, cylinder measures, and conditional cylinder functions. We further present a collection of technical lemmas and propositions to derive uniform quantitative bounds for Fourier frame analysis on overlapping self-similar measures. All results in this section are stated on $\R^d$.

\subsection{Basic Symbolic Dynamics and Attractor Geometry}
Let \(\{S_a\}_{a\in\D}\) be the IFS given in \eqref{eq:intro-ifs}. Write K for its attractor in \eqref{eq:intro-attractor} and \(\mu\) for the self-similar measure in \eqref{eq:intro-invariance} associated with this IFS. Let \(\Omega = \D^\mathbb{N}\) denote the full one-sided symbolic space over \(\D\) endowed with the product Bernoulli measure \eqref{eq:product-Bernoulli measure}.
Recall that the coding map \(\pi: \Omega\to K\) is defined by
\begin{equation*}
	\pi(\omega)
	:=\sum_{j=1}^{\infty}\rho^{j-1}b_{\omega_j}, \quad \omega=(\omega_j)_{j\ge 1}\in\Omega,
\end{equation*}
which is continuous and satisfies \(\mu=\pi_\#\Pp\). The identity \(\pi(\Omega)=\bigcup_{a\in\D}S_a(\pi(\Omega))\),  combined with the uniqueness of the IFS attractor, yields \(\pi(\Omega)=K\).  Consequently, \(\pi\) is surjective.  A
sequence $\omega\in\Omega$ with $\pi(\omega)=x$ is called a symbolic address
of $x\in K$.  Overlaps are precisely what allow one physical point to have
several addresses.  Since $\Pp$ has full support on $\Omega$ and $\pi$ is continuous, we deduce that $\supp\mu=K$.

\vspace{0.2cm}

Fix an integer $n\in\N$. For a finite word $u=(u_1,\ldots,u_n)\in\D^n$ and an infinite symbolic sequence $v=(v_1,v_2,\ldots)\in\Omega$, define the word length of $u$, the concatenation of $u$ and $v$, and the iterated affine map $S_u$, respectively, by
\begin{equation*}
  \abs u:=n,
  \qquad
  uv:=(u_1,\ldots,u_n,v_1,v_2,\ldots),
  \qquad
  S_u:=S_{u_1}\circ\cdots\circ S_{u_n}.
\end{equation*}
The translational shift associated with the word $u\in \D^n$ is
\[
  t_u:=\sum_{j=1}^n\rho^{j-1}b_{u_j}\in\R^d.
\]
Direct composition yields two critical compatibility relations for the affine maps and the coding map $\pi$:
\begin{equation}\label{eq:word-map-formula}
  S_u(x)=\rho^n x+t_u,
  \qquad
  \pi(u\tau)=S_u(\pi(\tau)).
\end{equation}

\subsection{Shannon Entropy and Level-\(n\) Fiber Entropy}

As $\Omega$ and $K$ are compact metric spaces, we apply Rokhlin's disintegration theorem (see, e.g., Kallenberg~\cite{Kallenberg2021}) to disintegrate the probability measure $\Pp$ with respect to the coding map $\pi$. This procedure yields a family of probability kernels 
$\{\Pp_x:x\in K\}$ with the following properties: each measure
$\Pp_x$ is a Borel probability measure on $\Omega$; the map
$x\mapsto\Pp_x(G)$ is Borel measurable for every $G\in\B(\Omega)$, and the integral representation
\begin{equation}\label{eq:intro-disintegration}
	\Pp(G\cap\pi^{-1}(E))=\int_E\Pp_x(G)\dd\mu(x)
\end{equation}
is valid for all \(G\in\B(\Omega)\) and \(E\in\B(K)\).
We further have 
\begin{equation}\label{eq:fiber-support}
	\Pp_x(\pi^{-1}(\{x\}))=1,
	\quad  \mu-a.e. \ x\in K.
\end{equation}
In probabilistic terms, $\Pp_x$ corresponds to  the conditional probability distribution of the symbolic
address given the physical point $x=\pi(\omega)$.

\vspace{0.2cm}

For $u=(u_1,\ldots,u_n)\in\D^n$, define its cylinder set by
\begin{equation*}
  [u]:=\{\omega\in\Omega:\omega_j=u_j\text{ for }1\le j\le n\}.
\end{equation*}
Then \(\Pp([u])=N^{-n}\). For \(x\in K\) and \(u\in\D^n\), define
\begin{equation*}
  p_u(x):=\Pp_x([u]),
  \qquad
  \cP_n:=\{[u]:u\in\D^n\}.
\end{equation*}
Thus $p_u(x)$ can be interpreted as the conditional probability that a symbolic address of $x$ begins with $u$, while $\cP_n$ records the level-$n$ cylinder partition of \(\Omega\). We refer to  \(p_u(x)\) as {\it{a conditional cylinder function}}.

\vspace{0.2cm}  

Adopt the convention \(0\log0=0\). For each \(x\in K\), {\it{the Shannon entropy of $\cP_n$}} with respect to $\Pp_x$ (see e.g.\cite{Walters1982}, Section 4.2) is defined by
\begin{equation}
	H_{\Pp_x}(\cP_n) = -\sum_{\abs u=n}p_u(x)\log p_u(x).
\end{equation}
We define the {\it {level-\(n\) fiber entropy}} as the average of this Shannon entropy with respect to the self-similar measure \(\mu\):
\begin{equation}\label{eq:Rn-fiber-entropy}
	\begin{aligned}
		R_n
		=\int_K H_{\Pp_x}(\cP_n)\dd\mu(x).
	\end{aligned}
\end{equation}

\subsection{Cylinder Measures}

For any finite word $w\in\D^n$, we define {\it{the cylinder measure \(\mu_w\)}} by
\begin{equation}\label{eq:cylinder-submeasure}
  \mu_w(E):=\Pp([w]\cap\pi^{-1}(E)),
  \qquad \forall E\in\B(K).
\end{equation}
The measure \(\mu_w\) denotes the total mass of \(\mu\) restricted to symbolic sequences with prefix \(w\). The level-$n$ cylinders partition $\Omega$, so
\begin{equation*}
  \mu=\sum_{\abs w=n}\mu_w,
  \qquad \mu_w\ll\mu,
\end{equation*}
where $\mu_w\ll\mu$ denotes the absolute continuity of $\mu_w$ with respect to $\mu$.

\vspace{0.2cm}

Recall that the conditional cylinder function is defined by
\begin{equation*}
  p_w(x):=\Pp_x([w]),
  \qquad x\in K.
\end{equation*}
The disintegration identity \eqref{eq:intro-disintegration} gives, for every $E\in\B(K)$,
\begin{equation}\label{eq:posterior-rn}
  \mu_w(E)=\int_Ep_w(x)\dd\mu(x),
  \qquad
  p_w=\frac{\dd\mu_w}{\dd\mu}\quad\mu\text{-a.e.}
\end{equation}
 $ p_w$ is therefore exactly the Radon--Nikodym density of the cylinder measure \(\mu_w\) with respect to \(\mu\). Since each $\Pp_x$ is a probability measure and $\cP_n$ forms a partition of $\Omega$, these functions satisfy
\begin{equation*}
  0\le p_w(x)\le1,
  \qquad
  \sum_{\abs w=n}p_w(x)=1,
  \qquad x\in K.
\end{equation*}
In particular, each $p_w$ belongs to $L^2(\mu)$.

\vspace{0.2cm} 

The following lemma provides an explicit representation of \(\mu_w\) and establishes a critical uniform Fourier magnitude identity that underlies all subsequent comparison arguments.

\begin{lemma}\label{lem:posterior}
For every $n\in\N$ and $w\in\D^n$,
\begin{equation}\label{eq:posterior-measure}
  \mu_w=N^{-n}(S_w)_\#\mu,
  \qquad
  \int_Kp_w\dd\mu=N^{-n},
\end{equation}
The conditional cylinder function \(p_w\) vanishes \(\mu\)-a.e. outside \(S_w(K)\), and its Fourier transform satisfies
\begin{equation}\label{eq:posterior-fourier}
  \wh{p_w\mu}(\xi)
  =N^{-n}e^{-2\pi i\xi\cdot t_w}\wh\mu(\rho^n\xi),
  \qquad
  \abs{\wh{p_w\mu}(\xi)}
  =N^{-n}\abs{\wh\mu(\rho^n\xi)},
  \quad\xi\in\R^d.
\end{equation}
\end{lemma}

\begin{proof}
 For every $E\in\B(K)$,	\eqref{eq:cylinder-submeasure} and the conditional probability formula gives
 \begin{equation}\label{eq:conditional-probability}
 \mu_w(E)=\Pp\big([w] \cap \pi^{-1}(E)\big) = \Pp([w]) \cdot \Pp\big(\pi(u)\in E \,\big|\, u\in [w]\big).
 \end{equation}
If \(u\in[w]\), we write \(u=wv\) for a unique infinite tail sequence \(v\in\Omega\). Under the Bernoulli product measure, conditioning on the prefix \(w\), the tail sequence \(v\) again has distribution \(\Pp\). The coding identity \(\pi(wv)=S_w(\pi(v))\) follows from \eqref{eq:word-map-formula}. We thereby compute
\begin{equation}\label{eq:word-map-probability}
 \Pp\big(\pi(u)\in E \,\big|\, u\in [w]\big)
 = \Pp\big(S_w(\pi(v)) \in E\big)
 = \Pp\big(\pi(v) \in S_w^{-1}(E)\big)
 = \mu\big(S_w^{-1}(E)\big).
\end{equation}
Combining \eqref{eq:conditional-probability} and \eqref{eq:word-map-probability} gives
\[
\mu_w(E)=N^{-n}\mu\big(S_w^{-1}(E)\big)
=N^{-n}\bigl((S_w)_\#\mu\bigr)(E).
\]
Integrating over the full attractor $K$, we obtain
\[
\int_Kp_w\dd\mu
=\mu_w(K)
=N^{-n}\bigl((S_w)_\#\mu\bigr)(K)
=N^{-n}\mu(K)
=N^{-n}.
\]
 Since $(S_w)_\#\mu$ is supported on $S_w(K)$,
\eqref{eq:posterior-rn} gives the support assertion for $p_w$.

 Finally, using $p_w\dd\mu=\dd\mu_w$ and the pushforward formula \eqref{eq:posterior-measure},
\[
\begin{aligned}
  \wh{p_w\mu}(\xi)
  &=\int_K e^{-2\pi i\xi\cdot y}\dd\mu_w(y)\\
  &=N^{-n}\int_K e^{-2\pi i\xi\cdot y}\dd\bigl((S_w)_\#\mu\bigr)(y)\\
  &=N^{-n}\int_K e^{-2\pi i\xi\cdot(\rho^n x+t_w)}\dd\mu(x)\\
  &=N^{-n}e^{-2\pi i\xi\cdot t_w}\wh\mu(\rho^n\xi).
\end{aligned}
\]
The proof is complete.
\end{proof}

\begin{remark}
	It is important to distinguish the cylinder measure $p_w\dd\mu$ from the geometric restriction $1_{K_w}\dd\mu$ of the set $K_w=S_w(K)$.
	Recall that $p_w$ denotes the Radon--Nikodym density of the symbolic cylinder measure $\mu_w$ with respect to $\mu$, and
	\begin{equation}\label{eq:remark-mu-w}
		p_w\dd\mu=\dd\mu_w=N^{-n}\dd\bigl((S_w)_\#\mu\bigr).
	\end{equation}
	For each $v\in\D^n$, $\mu_v:=N^{-n}(S_v)_\#\mu$ by \eqref{eq:remark-mu-w}.
	Since $\mu=\sum_{v\in\D^n}\mu_v$ and $\mu_w$ is supported on $K_w$, restricting $\mu$ to $K_w$ yields
		$\mu|_{K_w}
		=\mu_w+
		\sum_{\substack{v\in\D^n\\v\ne w}}\mu_v|_{K_w}$.
	Consequently,
		\[
		p_w=1_{K_w}\quad\mu\text{-a.e.}
		\quad\Longleftrightarrow\quad
		\mu_v(K_w)=0\quad\text{for every }v\in\D^n\setminus\{w\}.
		\]
	In particular, $p_w\dd\mu=1_{K_w}\dd\mu$ whenever the level-$n$ geometric cylinders are pairwise disjoint modulo $\mu$-null sets.
	Nonetheless, a nonempty geometric intersection by itself does not force the two measures to differ.
	If another branch carries positive measure into $K_w$, then $\mu|_{K_w}\ne\mu_w$.
	In fact, their Fourier transforms already differ at the origin:
		\[
		\begin{aligned}
			\widehat{1_{K_w}\mu}(0)-\widehat{p_w\mu}(0)
			&=\mu(K_w)-N^{-n}\\
			&=\sum_{\substack{v\in\D^n\\v\ne w}}\mu_v(K_w)>0.
		\end{aligned}
		\]
	In contrast, the pushforward identity \eqref{eq:remark-mu-w} always gives
	$\widehat{p_w\mu}(\xi)
	=N^{-n}e^{-2\pi i\xi\cdot t_w}\widehat\mu(\rho^n\xi)$.
	The modulus of the Fourier transform of $p_w \dd\mu$ is therefore independent of the word $w$ at a fixed level, while the corresponding modulus identity generally fails for $\widehat{1_{K_w}\mu}$.
	This distinction is exactly the reason we use conditional cylinder densities, rather than geometric set indicators, in our Fourier frame analysis.
\end{remark}

Throughout the remainder of the paper, we assume that $\Gamma$ is a frame measure for $\mu$ with lower and upper frame bounds
$A,B$, respectively.  The modulus of the right-hand side of
\eqref{eq:posterior-fourier} is independent of $u$ at each fixed level.
Applying \eqref{eq:frame-inequality} to two such conditional cylinder functions gives
the following comparison.

\begin{proposition}\label{prop:comparability}
Assume that \(\Gamma\) is a frame measure for \(\mu\) with lower and upper frame bounds \(A\) and \(B\). Then, for every integer \(n\ge1\) and every pair \(v,u\in\D^n\),
\begin{equation}\label{eq:comparability}
  \frac{A}{B}\norm{p_v}_{L^2(\mu)}^2
  \le \norm{p_u}_{L^2(\mu)}^2
  \le \frac{B}{A}\norm{p_v}_{L^2(\mu)}^2.
\end{equation}
\end{proposition}

\begin{proof}
By \eqref{eq:posterior-fourier}, the quantity
\begin{equation*}
  \Theta_n
  :=\int_{\R^d}\abs{\wh{p_u\mu}(\xi)}^2\dd\Gamma(\xi)
  =N^{-2n}\int_{\R^d}\abs{\wh\mu(\rho^n\xi)}^2\dd\Gamma(\xi)
\end{equation*}
is independent of the word $u$ with $\abs u=n$.  For either \(z=v\) or \(z=u\), the frame inequalities give
\begin{equation*}
  A\norm{p_z}_{L^2(\mu)}^2
  \le\Theta_n
  \le B\norm{p_z}_{L^2(\mu)}^2.
\end{equation*}
Combining the lower bound for $z=v$ with the upper bound for $z=u$ yields
\begin{equation*}
  \frac{A}{B}\norm{p_v}_{L^2(\mu)}^2
  \le\norm{p_u}_{L^2(\mu)}^2,
\end{equation*}
and reversing the roles of $v$ and $u$ yields the other inequality in
\eqref{eq:comparability}.
\end{proof}

\section{Collision Bounds and Atomicity}\label{sec:atomicity}
We first use the rightmost branch in dimension one to establish a uniform lower bound for conditional collisions. We then show that such a bound forces the conditional coding measures to be purely atomic almost everywhere. The second argument applies in arbitrary dimension and will also be used in Section~\ref{sec:Rd}.

\subsection{The rightmost branch and the collision lower bound}\label{subsec:collision-estimate}

In this subsection, we specialize to the case $d=1$. After relabeling the
translation parameters, assume that
\[
  b_0<b_1<\cdots<b_{N-1}.
\]
Set
\[
  x_-:=\min K=\frac{b_0}{1-\rho},
  \qquad
  x_+:=\max K=\frac{b_{N-1}}{1-\rho}.
\]
Indeed, these extremal values are obtained by the two constant infinite symbolic sequences and 
$x_-\le \pi(\omega)\le x_+$ for every $\omega\in\Omega$.
For $n\in\N$, let 
\[ u_n^+ := (\underbrace{N-1,N-1,\dots,N-1}_{n\text{ entries}})\in\D^n,
\qquad I_n := S_{u_n^+}(K) = S_{N-1}^{\circ n}(K),\]
where $S_{N-1}^{\circ n}$ denotes the $n$-fold composition of $S_{N-1}$. Lemma \ref{lem:endpoint} below has been proved in~\cite[Proposition~3.9]{DutkayLai2014}, we include it here for later use.

\begin{lemma}[Proposition~3.9,\cite{DutkayLai2014}]\label{lem:endpoint}
There exist $n_0\ge1$ and $C_0>0$ such that for all integers \(n\ge n_0\),
\begin{equation*}\label{eq:endpoint-mass}
  \mu(I_n)=C_0N^{-n}.
\end{equation*}
\end{lemma}

For $n\in\N$ and $x\in K$, define 
\begin{equation*}
  c_n(x):=\sum_{\abs u=n}(\Pp_x([u]))^2
  =\sum_{\abs u=n}(p_u(x))^2.
\end{equation*}
An equivalent probabilistic interpretation is as follows: if \(\omega,\omega'\) are independent random symbolic sequences sampled according to \(\Pp_x\), then 
\begin{equation*}
  c_n(x)
  =(\Pp_x\otimes\Pp_x)
    \{(\omega,\omega'):\omega_j=\omega'_j\text{ for }1\le j\le n\}.
\end{equation*}
The quantity \(c_n(x)\) satisfies \(0\le c_n(x)\le 1\) pointwise. Combining Lemma \ref{lem:endpoint} and Proposition \ref{prop:comparability}, we derive a uniform positive lower bound for the integral of \(c_n\) over $K$.  

\begin{proposition}\label{prop:collision-lower}
There is $\kappa>0$ such that
\begin{equation}\label{eq:mean-collision-lower}
  \int_Kc_n(x)\dd\mu(x)\ge\kappa,
  \qquad \forall n\ge n_0.
\end{equation}
\end{proposition}

\begin{proof}
It follows from Lemma \ref{lem:posterior} that the conditional cylinder function $p_{u_n^+}$ is supported on $I_n$ with total integral \(\int_K p_{u_n^+}\,d\mu = N^{-n}\).
Hence Cauchy--Schwarz  inequality gives
\[
\begin{aligned}
  N^{-n}
  &=\int_{I_n}p_{u_n^+}(x)\dd\mu(x)\\
  &\le (\mu(I_n))^{1/2}
     \left(\int_{I_n}(p_{u_n^+}(x))^2\dd\mu(x)\right)^{1/2}\\
  &\le (\mu(I_n))^{1/2}\norm{p_{u_n^+}}_{L^2(\mu)}.
\end{aligned}
\]
After squaring both sides and applying Lemma~\ref{lem:endpoint}, we obtain
\[
  N^{-2n}
  \le C_0N^{-n}\norm{p_{u_n^+}}_{L^2(\mu)}^2,
\]
so $\norm{p_{u_n^+}}_{L^2(\mu)}^2\ge C_0^{-1}N^{-n}$.  By
Proposition~\ref{prop:comparability}, every length-$n$ conditional cylinder function satisfies
\[
  \norm{p_u}_{L^2(\mu)}^2\ge\frac{A}{BC_0}N^{-n},
  \qquad \forall \ \abs u=n.
\]
Summing over all \(N^n\) distinct length-$n$ words gives
\[
\begin{aligned}
  \int_K c_n(x)\dd\mu(x)
  &=\sum_{\abs u=n}\norm{p_u}_{L^2(\mu)}^2\\
  &\ge N^n\frac{A}{BC_0}N^{-n}
   =\frac{A}{BC_0}.
\end{aligned}
\]
Set $\kappa=A/(BC_0)$. This completes the proof of \eqref{eq:mean-collision-lower}.
\end{proof}

\subsection{Pure atomicity of almost every conditional coding measure}

Proposition~\ref{prop:collision-lower} verifies
\eqref{eq:mean-collision-lower} in dimension one, while the
higher-dimensional verification will be supplied in
Section~\ref{sec:Rd}. The remaining argument applies in arbitrary dimension. Assume that \eqref{eq:mean-collision-lower} holds in $\R^d$.
For any probability measure \(\eta\) on the symbolic space \(\Omega\), define its atomic support set by
\begin{equation*}
  \At(\eta):=\{w\in\Omega:\eta(\{w\})>0\}.
\end{equation*}
The set $\At(\eta)$ is at most countable.  A probability measure $\eta$ is purely
atomic if and only if
\begin{equation*}
  \eta(\At(\eta))
  =\sum_{w\in\At(\eta)}\eta(\{w\})=1.
\end{equation*}

\vspace{0.2cm}

The following lemma links the quantity $c_n(x)$ to the sum of squared atomic masses associated with each conditional coding measure \(\Pp_x\).

\begin{lemma}\label{lem:collision-limit}
For every $x\in K$, the sequence  $c_n(x)$ is nonincreasing, and
\begin{equation*}
  \lim_{n\to\infty}c_n(x)
  =\sum_{w\in\At(\Pp_x)}(\Pp_x(\{w\}))^2.
\end{equation*}
In addition,
\begin{equation}\label{eq:positive-collision-limit}
  \int_K\lim_{n\to\infty}c_n(x)\dd\mu(x)\ge\kappa.
\end{equation}
\end{lemma}

\begin{proof}
Any length-$n$ cylinder \([u]\) decomposes as a disjoint union of length-\((n+1)\) cylinders: 
\[[u]=\bigsqcup_{a\in\D}[ua].\] 
By the convexity of the square function,
\[
  \sum_{a\in\D}(\Pp_x([ua]))^2
  \le\left(\sum_{a\in\D}\Pp_x([ua])\right)^2
  =(\Pp_x([u]))^2.
\]
Summing over all words $u$ of length $n$ immediately gives $c_{n+1}(x)\le c_n(x)$, so \(\{c_n(x)\}\) decreases pointwise in $n$.  
Now define the product sets
\[
  D_n:=\bigsqcup_{\abs u=n}[u]\times[u],
  \qquad
  \Delta_\Omega:=\{(w,w):w\in\Omega\},
\]
where \(\Delta_\Omega\) denotes the diagonal of the symbolic product space \(\Omega\times\Omega\).
The sets \(D_n\) form a decreasing sequence converging downwards to the diagonal: $D_n\downarrow\Delta_\Omega$. By the definition of \(c_n(x)\), we have 
\[(\Pp_x\otimes\Pp_x)(D_n)=c_n(x).\]  
Applying continuity from above to the product measure \(\Pp_x\otimes\Pp_x\), we obtain
\[
  \lim_{n\to\infty}c_n(x)
  =(\Pp_x\otimes\Pp_x)(\Delta_\Omega)
  =\sum_{w\in\At(\Pp_x)}(\Pp_x(\{w\}))^2.
\]
Finally, the Bounded Convergence Theorem together with 
\eqref{eq:mean-collision-lower} yields
\eqref{eq:positive-collision-limit}.
\end{proof}

The inequality \eqref{eq:positive-collision-limit} in Lemma \ref{lem:collision-limit} corresponds to a measurable set of symbolic addresses that carry strictly positive point masses within their respective conditional coding measures.

\begin{lemma}\label{lem:positive-atom-event}
There exists a Borel subset $\mathcal A\subset\Omega$ such that
\begin{equation}\label{eq:fiber-atomic-mass}
  \Pp_x(\mathcal A)
  =\sum_{w\in\At(\Pp_x)}\Pp_x(\{w\})
\end{equation}
for $\mu$-almost every $x\in K$, with the global measure bound $\Pp(\mathcal A)>0$.
\end{lemma}

\begin{proof}
For each integer $n\ge1$, define 
\begin{equation*}
  \alpha_n(w)
  :=\sum_{\abs u=n}\1_{[u]}(w)p_u(\pi(w))
  =\Pp_{\pi(w)}([w_1\cdots w_n]).
\end{equation*}
The nested decreasing cylinder sets \([w_1\cdots w_n]\) converges to the singleton \(\{w\}\), so we can take the pointwise limit
\begin{equation*}
  \alpha(w):=\lim_{n\to\infty}\alpha_n(w)
  =\Pp_{\pi(w)}(\{w\}).
\end{equation*}
Each \(\alpha_n\) is Borel measurable, so \(\alpha\) is also Borel measurable. Define
\[\mathcal{A} := \{w\in\Omega: \alpha(w)>0\}.\]
By the support
property of $\Pp_x$ in \eqref{eq:fiber-support}, for $\mu$-almost every $x\in K$, the set
$\mathcal A$ coincides with  $\At(\Pp_x)$ up to a $\Pp_x$-null set. This immediately yields equation
\eqref{eq:fiber-atomic-mass}.  Since $t^2\le t$ for every $t\in[0,1]$, applying this together with Lemma~\ref{lem:collision-limit} and the disintegration formula \(\Pp=\int_K \Pp_x\,d\mu(x)\), we compute
\[\Pp(\mathcal{A}) = \int_K \Pp_x(\mathcal{A})\,d\mu(x)
\ge \int_K \sum_{w\in\operatorname{At}(\Pp_x)} (\Pp_x(\{w\}))^2 d\mu(x)
= \int_K \lim_{n\to\infty} c_n(x)\,d\mu(x) \ge \kappa.\]
Since \(\kappa>0\), we conclude \(\Pp(\mathcal{A})>0\).
\end{proof}

Fix a point $x^\circ\in K$. For each digit
$a\in\D$, we define $r_a:K\to[0,1]$ and
$T_a:K\to K$ by
\begin{equation*}
  r_a(x):=\Pp_x([a]),
  \qquad
  T_a(x):=
  \begin{cases}
    S_a^{-1}(x),&x\in S_a(K),\\
    x^\circ,&x\notin S_a(K).
  \end{cases}
\end{equation*}
Both \(r_a\) and \(T_a\) are Borel measurable. 
Let \(\sigma\colon\Omega\to\Omega\) denote the left-shift operator on the symbolic space defined by
\[\sigma(\omega_1,\omega_2,\omega_3,\dots)=(\omega_2,\omega_3,\omega_4,\dots).\]
In order to push forward positive atomic masses along the Bernoulli shift dynamics, we establish a Bayes-style factorization formula which holds for all Borel events $G$ outside a $\mu$-null exceptional set (this null set may be chosen independently of $G$).

\begin{proposition}
\label{prop:bayes}
For each $a\in\D$, there exists a Borel measurable set
$\mathcal N_a\subset K$ with $\mu(\mathcal N_a)=0$ such that, for every
$x\in K\setminus\mathcal N_a$ and every $G\in\B(\Omega)$,
\begin{equation}\label{eq:bayes}
  \Pp_x([a]\cap\sigma^{-1}(G))
  =\1_{S_a(K)}(x)r_a(x)\Pp_{T_a(x)}(G).
\end{equation}
\end{proposition}

\begin{proof}
Fix \(a\in\D\) and a Borel set \(G\in\B(\Omega)\). For $E\in\B(K)$, combining the formula \eqref{eq:intro-disintegration} and the identity \(\pi(a\tau)=S_a(\pi(\tau))\), we have 
\begin{equation*}
\begin{aligned}
  \int_E\Pp_x([a]\cap\sigma^{-1}(G))\dd\mu(x)
  &=\frac1N\int_\Omega
      \1_G(\tau)\1_E(S_a(\pi(\tau)))\dd\Pp(\tau)\\
  &=\frac1N\int_K
      \1_E(S_a(y))\Pp_y(G)\dd\mu(y).
\end{aligned}
\end{equation*}
Recall that $\mu_a(E)=\Pp([a]\cap\pi^{-1}(E))$.  By
Lemma~\ref{lem:posterior},
$\mu_a=N^{-1}(S_a)_\#\mu$ and $r_a=\dd\mu_a/\dd\mu$.  Applying
$r_a\dd\mu=\dd\mu_a$ to the Borel function
$\1_E\1_{S_a(K)}\Pp_{T_a(\cdot)}(G)$ gives
\begin{equation*}
\begin{aligned}
  &\int_E\1_{S_a(K)}(x)r_a(x)\Pp_{T_a(x)}(G)\dd\mu(x)\\
  &\qquad=\frac1N\int_K
      \1_E(S_a(y))\Pp_{T_a(S_a(y))}(G)\dd\mu(y)\\
  &\qquad=\frac1N\int_K
      \1_E(S_a(y))\Pp_y(G)\dd\mu(y).
\end{aligned}
\end{equation*}
Both sides of \eqref{eq:bayes} are nonnegative and their integrals over arbitrary Borel sets \(E\) coincide, so \eqref{eq:bayes} holds for \(\mu\)-a.e.  $x\in K$ with a null exceptional set
 depending only on the fixed $G$.  Let \(\mathcal{C}_0\) denote the countable algebra generated by all finite unions of symbolic cylinder sets, and let \(\mathcal{N}_{a,G}\subset K\) be the null set where \eqref{eq:bayes} fails for a  fixed \(G\in\mathcal{C}_0\).  Set
\(
  \mathcal N_a:=\bigcup_{G\in\mathcal C_0}\mathcal N_{a,G}\). Then
\(
  \mu(\mathcal N_a)=0
\)
as a countable union of \(\mu\)-null sets remains null. 
For any fixed $x\notin\mathcal N_a$, both sides of \eqref{eq:bayes} define finite Borel measures in the variable \(G\), and these two measures agree on the generating algebra \(\mathcal{C}_0\). By the monotone class theorem, the equality extends to every Borel set \(G\in\B(\Omega)\).
\end{proof}

We now combine the positive-mass atomic event constructed in Lemma \ref{lem:positive-atom-event} above, the recursive Bayes factorization from Proposition \ref{prop:bayes}, and ergodicity of the Bernoulli shift to deduce pure atomicity for nearly all conditional coding measures.

\begin{proposition}\label{prop:pure-atomic}
For $\mu$-almost every $x\in K$, the conditional coding measure $\Pp_x$ is
purely atomic.
\end{proposition}

\begin{proof}
We first verify that the conditional probability of the leading digit is strictly positive for almost every symbolic sequence. For each digit \(a\in\D\), define the zero level set
 \[\mathcal{Z}_a=\{x:r_a(x)=0\}.\]
By the disintegration formula \eqref{eq:intro-disintegration},
\[
  \Pp\!\left(\bigcup_{a\in\D}([a]\cap\pi^{-1}(\mathcal{Z}_a))\right)
  =\sum_{a\in\D}\int_{\mathcal{Z}_a}r_a(x)\dd\mu(x)=0.
\]
This immediately yields
\begin{equation*}
  r_{w_1}(\pi(w))>0
  \quad\text{for $\Pp$-almost every } w=(w_1,w_2,\ldots)\in\Omega.
\end{equation*}

Let $\mathcal N=\bigcup_{a\in\D}\mathcal N_a$, where \(\mathcal N_a\) is defined  in Proposition \ref{prop:bayes}. Then
$\mu(\mathcal N)=0$, so $\pi(w)\notin\mathcal N$ for almost every
$w$.  Write $w=w_1w'$ with the leading digit  $w_1$ and the shifted tail 
$w'=\sigma(w)$. The coding map satisfies two fundamental identities:
\[
  \pi(w)=S_{w_1}(\pi(w')),
  \qquad
  \{w\}=[w_1]\cap\sigma^{-1}(\{w'\}).
\]
Applying Proposition~\ref{prop:bayes} to the singleton set \(G=\{w'\}\), for \(\Pp\)-almost every \(w\), we have
\begin{equation}\label{eq:alpha-recursion}
\begin{aligned}
  \alpha(w)
  &=\Pp_{\pi(w)}(\{w\})\\
  &=\Pp_{\pi(w)}\bigl([w_1]\cap\sigma^{-1}\{w'\}\bigr)\\
  &=r_{w_1}(\pi(w))\Pp_{T_{w_1}(\pi(w))}(\{w'\})\\
  &=r_{w_1}(\pi(w))\alpha(\sigma(w)),
\end{aligned}
\end{equation}
where $T_{w_1}(\pi(w))=\pi(w')$ by
\eqref{eq:word-map-formula}.
Since the factor $r_{w_1}(\pi(w))$ is positive almost surely,
\eqref{eq:alpha-recursion} gives the key equivalence
\begin{equation*}
  \alpha(w)>0
  \quad\Longleftrightarrow\quad
  \alpha(\sigma(w))>0
  \qquad\text{for $\Pp$-a.e. } w\in\Omega.
\end{equation*}
Equivalently, for  \(\mathcal{A}=\{w\in\Omega:\alpha(w)>0\}\),
\begin{equation*}
  \Pp\bigl(\mathcal A\mathbin{\triangle}\sigma^{-1}(\mathcal A)\bigr)=0.
\end{equation*}
Thus \(\mathcal A\) is shift-invariant modulo null sets. The Bernoulli shift is ergodic (see  \cite{Walters1982}), which forces every shift-invariant measurable set to have global measure either $0$ or $1$.  Lemma~\ref{lem:positive-atom-event} gives
$\Pp(\mathcal A)>0$, hence $\Pp(\mathcal A)=1$.  Finally, using \eqref{eq:intro-disintegration} again, we have
\[
  1=\Pp(\mathcal A)
  =\int_K\sum_{w\in\At(\Pp_x)}
      \Pp_x(\{w\})\dd\mu(x).
\]
This implies that $\sum_{w\in\At(\Pp_x)}
\Pp_x(\{w\}) =1$ 
$\mu$-a.e. $x\in K$ since it is bounded by $1$. Thus each \(\Pp_x\) is purely atomic for \(\mu\)-a.e. $x\in K$.
\end{proof}

\Needspace{10\baselineskip}
\section{Dirac Rigidity and the Dimension Formula}\label{sec:entropy}
This section completes the rigidity argument and derives its dimensional consequences. We first obtain the Dirac property and essential injectivity of the coding map, and then apply the Feng--Hu projection-entropy framework.

\subsection{Essential injectivity of the coding map}
\label{subsec:dirac-essential}

Proposition~\ref{prop:pure-atomic} establishes pure atomicity for generic conditional coding measures. We now derive a stronger result from the Bayes-type identity \eqref{eq:bayes}: nearly every conditional coding measure is supported on a single symbolic atom. 

\vspace{0.2cm}

Recall the Borel function
$
  \alpha(w)=\Pp_{\pi(w)}(\{w\})
$
and the positive-mass atomic event \(\mathcal{A} = \{w\in\Omega: \alpha(w)>0\}\)  introduced in
Lemma~\ref{lem:positive-atom-event}.  The proof of
Proposition~\ref{prop:pure-atomic} gives $\Pp(\mathcal A)=1$, and so
\begin{equation}\label{eq:alpha-positive}
  \alpha(w)>0
  \quad\text{for $\Pp$-almost every }w\in\Omega.
\end{equation}
Shift invariance of \(\Pp\) will further force the singleton mass functional \(\alpha\) to be stationary along \(\Pp\)-almost every orbit, as formalized in the following proposition.

\begin{proposition}
\label{prop:stationary-singleton}
For $\Pp$-almost every $w\in\Omega$,
\begin{equation}\label{eq:stationary-singleton}
  \alpha(w)=\alpha(\sigma(w)),
  \qquad
  r_{w_1}(\pi(w))=1.
\end{equation}
\end{proposition}

\begin{proof}
The identity \eqref{eq:alpha-recursion} and
$0\le r_{w_1}(\pi(w))\le1$ imply that
\[
  0\le\alpha(w)\le\alpha(\sigma(w))
  \qquad\text{for $\Pp$-a.e.} \ w\in\Omega.
\]
Since \(0\le \alpha (w)\le 1\), we have \(\alpha(w)\in L^1(\Pp)\). The Bernoulli measure \(\Pp\) is shift-invariant, so
\[
  \int_\Omega \alpha(\sigma(w))\dd\Pp(w)
  =\int_\Omega \alpha(w)\dd\Pp(w).
\]
This shows that the nonnegative function $\alpha\circ\sigma-\alpha$ has integral zero, and consequently
$\alpha=\alpha\circ\sigma$ $\Pp$-almost everywhere.  Substituting this equality
into \eqref{eq:alpha-recursion} gives
\[
  \alpha(w)
  =r_{w_1}(\pi(w))\alpha(w).
\]
By \eqref{eq:alpha-positive},  the second equality in
\eqref{eq:stationary-singleton} follows.
\end{proof}

The following lemma iterates the single-step Bayes factorization identity \eqref{eq:bayes} established in Proposition \ref{prop:bayes}. Distinct null sets for different cylinders cannot support a globally valid product formula for all \(n\ge1\). We therefore construct a full-measure shift-invariant set \(\Omega_\ast\) as a unified domain for the identity.

\begin{lemma}
\label{lem:iterated-posterior}
There exists a forward-invariant Borel set
$\Omega_\ast\subset\Omega$ with $\Pp(\Omega_\ast)=1$ such that, for every
$w\in\Omega_\ast$ and every $n\ge1$,
\begin{equation}\label{eq:iterated-posterior}
  \Pp_{\pi(w)}([w_1\cdots w_n])
  =\prod_{j=0}^{n-1}
    r_{w_{j+1}}\bigl(\pi(\sigma^j(w))\bigr)=1.
\end{equation}
\end{lemma}

\begin{proof}
Let
\(
  \mathcal N:=\bigcup_{a\in\D}\mathcal N_a,
\)
where each Borel null set $\mathcal N_a\subset K$ denotes the exceptional set obtained in
Proposition~\ref{prop:bayes}. Define 
$$
  E_{\mathrm B}:=\{w\in\Omega:\pi(w)\notin\mathcal N\}.
$$
Then \(E_{\mathrm B}\) is Borel measurable as a preimage of a Borel set under the measurable coding map \(\pi\). Recall that $\mu=\pi_\#\Pp$ and note that $\mu(\mathcal N)=0$. The pushforward property forces \(\Pp(E_{\mathrm B})=1\). Set
\[
  E_{\mathrm R}
  :=\{w\in\Omega:
      r_{w_1}(\pi(w))=1\}.
\]
Proposition~\ref{prop:stationary-singleton} guarantees \(E_{\mathrm R}\) is also a Borel set of  full \(\Pp\)-measure. Let
\[
  \Omega_\ast
  :=\bigcap_{j=0}^{\infty}
       \sigma^{-j}(E_{\mathrm B}\cap E_{\mathrm R}).
\]
Then $\Omega_\ast$ is Borel measurable, $\Pp(\Omega_\ast)=1$, and
$\sigma(\Omega_\ast)\subset\Omega_\ast$, i.e., \(\Omega_\ast\) is forward invariant under the shift map. Moreover, for every
$w\in\Omega_\ast$ and $j\ge0$,
$r_{w_{j+1}}\bigl(\pi(\sigma^j(w))\bigr)=1
$ by Proposition \ref{prop:stationary-singleton}.

We verify the product equality~\eqref{eq:iterated-posterior}  by induction on cylinder length \(n\).
For $n=1$, the left-hand side  \(\Pp_{\pi(w)}([w_1]) = r_{w_1}(\pi(w)) =1\). 
For $n\ge2$, set $G=[w_2\cdots w_n]$. The cylinder decomposition
\[
  [w_1\cdots w_n]
  =[w_1]\cap\sigma^{-1}(G),
\]
holds by definition of the shift operator on symbolic cylinders. Applying Proposition~\ref{prop:bayes} at $x=\pi(w)$ gives
\[
  \Pp_{\pi(w)}([w_1\cdots w_n])
  =r_{w_1}(\pi(w))
   \Pp_{\pi(\sigma(w))}([w_2\cdots w_n]).
\]
Since $\sigma^j(w)\in\Omega_\ast$ for every $j\ge0$, the inductive iteration is well-defined at every shift level. Every factor in the resulting finite product equals $1$ by construction of \(\Omega_\ast\), completing the proof of~\eqref{eq:iterated-posterior}.
\end{proof}

The iterated product identity derived in Lemma~\ref{lem:iterated-posterior} restricts the full mass of each conditional coding measure to exactly one symbolic address.

\begin{proposition}
\label{prop:dirac-fibers}
For $\mu$-almost every $x\in K$, the conditional coding measure $\Pp_x$ is
a Dirac mass supported on a single symbolic sequence.
\end{proposition}

\begin{proof}
Fix any \(w\in\Omega_\ast\), where $\Omega_\ast$ is defined in Lemma~\ref{lem:iterated-posterior}. The nested cylinder sets \([w_1\cdots w_n]\) form a decreasing sequence converging to the singleton \(\{w\}\) as \(n\to\infty\). By the continuity from above property of finite measures and Lemma~\ref{lem:iterated-posterior}, 
\begin{equation}\label{eq:singleton-one}
  \Pp_{\pi(w)}(\{w\})
  =\lim_{n\to\infty}
    \Pp_{\pi(w)}([w_1\cdots w_n])=1.
\end{equation}
Disintegrating the full-\(\Pp\)-measure set \(\Omega_\ast\) via the coding map decomposition  gives
\[
  1=\Pp(\Omega_\ast)
  =\int_K\Pp_x(\Omega_\ast)\dd\mu(x).
\]
Then \(\Pp_x(\Omega_\ast)=1\) for \(\mu\)-a.e. \(x\in K\) since \(0\le \Pp_x(\Omega_\ast)\le 1\) pointwise for all $x\in K$. Denote
$$ K_\ast = \{x\in K: \Pp_x(\Omega_\ast)=1 \ \text{and} \ \Pp_x(\pi^{-1}(\{x\}))=1\}.$$
Then $\mu(K_\ast)=1$ by \eqref{eq:fiber-support} and the above argument.  For such $x\in K_\ast$, pick a symbolic sequence 
$w\in\Omega_\ast\cap\pi^{-1}(\{x\})$, which exists since \(\Pp_x(\Omega_\ast)=1\). \eqref{eq:singleton-one} then implies
$\Pp_x(\{\omega\})=1$, and so the conditional coding measure satisfies $\Pp_x=\delta_w$.
\end{proof}


\begin{corollary}\label{cor:exact-entropy}
For every $n\ge1$, \( R_n=0\).
\end{corollary}

\begin{proof}
 $H_{\Pp_x}(\cP_n)=0$ almost everywhere by
Proposition~\ref{prop:dirac-fibers}. Integration in
\eqref{eq:Rn-fiber-entropy} gives the statement.
\end{proof}

The Dirac structure of conditional coding measures further allows us to formalize the essential injectivity of the coding projection \(\pi\).

\begin{corollary}
\label{cor:essential-injectivity}
There exist Borel measurable subsets $\Omega_0\subset\Omega$ and $K_0\subset K$ satisfying
\[
  \Pp(\Omega_0)=\mu(K_0)=1
\]
such that the restricted coding map
\[
  \pi:\Omega_0\longrightarrow K_0
\]
is a Borel bijection.
\end{corollary}

\begin{proof}
Recall the Borel functional
\(
  \alpha(w)=\Pp_{\pi(w)}(\{w\}).
\)
 Equation~\eqref{eq:singleton-one} shows that the full-measure
set $\Omega_\ast$ is contained in
\[
  \Omega_0:=\{w\in\Omega:\alpha(w)=1\},
\]
so $\Pp(\Omega_0)=1$.  If $w, u\in\Omega_0$ and
$\pi(w)=\pi(u)=x$, then by the definition of \(\Omega_0\),
\[
  \Pp_x(\{w\})=\Pp_x(\{u\})=1.
\]
A probability measure can assign full unit mass to at most one singleton set, so we necessarily obtain $w=u$.  This proves the restriction $\pi|_{\Omega_0}$ is injective.

The set \(\Omega_0\) is Borel by construction, $\Omega$ and $K$ are Polish spaces, and the coding map $\pi$
is continuous, hence Borel measurable. Applying the Lusin--Souslin theorem,  the image
\[
  K_0:=\pi(\Omega_0)
\]
is a Borel subset of $K$, and the inverse mapping
\[
  \vartheta:=(\pi|_{\Omega_0})^{-1}:K_0\longrightarrow\Omega_0
\]
is also Borel measurable.  Finally, since \(K_0=\pi(\Omega_0)\) and \(\Omega_0\subset\pi^{-1}(K_0)\), the pushforward relation \(\mu=\pi_\#\Pp\) gives
\[
  \mu(K_0)=\Pp\bigl(\pi^{-1}(K_0)\bigr)\ge \Pp(\Omega_0)=1.
\]
Hence \(\mu(K_0)=1\), completing the proof.
\end{proof}

\begin{remark}
 This Borel inverse map \(\vartheta\) encodes the measure-theoretic disjointness of equal-depth symbolic cylinders. 
For any finite word \(u\in\D^n\), define
\[
E_u:=K_0\cap\vartheta^{-1}([u]).
\]
The collection \(\{E_u\}_{|u|=n}\) forms a Borel partition of $K$ up to a \(\mu\)-null set. Fix \(x\in E_u\), and write \(\vartheta(x)=uw\) for some tail sequence \(w\in\Omega\). Applying the fundamental coding relation yields
\[
x=\pi(\vartheta(x))=\pi(uw)=S_u(\pi(w))\in S_u(K),
\]
which implies that \(E_u \subset S_u(K)\). Furthermore, for every \(x\in K_0\), the definition of \(\Omega_0\) gives \(\Pp_x = \delta_{\vartheta(x)}\). It then follows that
\[p_u(x) = \Pp_x([u]) = \mathbf{1}_{E_u}(x).\]
Combining this equality with \eqref{eq:posterior-rn} and
Lemma~\ref{lem:posterior}, we obtain the restriction formula
\begin{equation}\label{eq:exact-cylinder-measure}
	\mu|_{E_u}=p_u  d\mu=\mu_u=N^{-n}(S_u)_\#\mu,
	\quad
	\mu(E_u)=N^{-n}.
\end{equation}
In particular, the pushforward measures from distinct equal-length words are mutually singular, i.e.,
\[(S_v)_\#\mu \perp (S_u)_\#\mu \quad \text{whenever } v\neq u \text{ and } |v|=|u|.\]
\end{remark}

\subsection{Projection Entropy and the Dimension Formula}

We now use the dimension theorem of Feng and Hu~\cite[Theorem~2.8]{FengHu2009} to obtain the remaining conclusion of Theorem~\ref{thm:rigidity}. Only the first-level conditional entropy is needed.

Let
\[
  \G:=\pi^{-1}(\B(K)),\qquad
  \cP:=\{[a]:a\in\D\}.
\]
For a finite Borel partition $\xi$ and a sub-$\sigma$-algebra
$\mathcal H$, denote the conditional Shannon entropy by
\[
  H_\Pp(\xi\mid\mathcal H)
  :=-\int_\Omega\sum_{C\in\xi}q_C\log q_C\,\dd\Pp,
  \qquad q_C:=\E_\Pp(\1_C\mid\mathcal H),
\]
with $0\log0=0$. Recall that the projection entropy of Feng and Hu is
\begin{equation}\label{eq:projection-entropy-definition}
  h_\pi(\sigma,\Pp)
  :=H_\Pp(\cP\mid\sigma^{-1}\G)
     -H_\Pp(\cP\mid\G).
\end{equation}

\begin{proposition}\label{prop:feng-hu-interface}
The measure $\mu$ is exact dimensional, and
\begin{equation}\label{eq:feng-hu-specialized}
  \dimH\mu=\frac{\log N-R_1}{-\log\rho}.
\end{equation}
\end{proposition}

\begin{proof}
For every $a\in\D$ and $E\in\B(K)$, the disintegration identity
\eqref{eq:intro-disintegration} gives
\[
\int_{\pi^{-1}(E)}\1_{[a]}\dd\Pp
  =\int_E p_a(x)\dd\mu(x)
  =\int_{\pi^{-1}(E)}p_a(\pi(\omega))\dd\Pp(\omega).
\]
Since $p_a\circ\pi$ is $\G$-measurable, this proves
$\E_\Pp(\1_{[a]}\mid\G)=p_a\circ\pi$ almost everywhere. Consequently,
\[
  H_\Pp(\cP\mid\G)
  =-\int_K\sum_{a\in\D}p_a(x)\log p_a(x)\dd\mu(x)
  =R_1.
\]
On the other hand, $\sigma^{-1}\G$ depends only on the tail
$(\omega_2,\omega_3,\ldots)$, which is independent of $\omega_1$ under
the product Bernoulli measure $\Pp$. Therefore
\[
  H_\Pp(\cP\mid\sigma^{-1}\G)
  =H_\Pp(\cP)=\log N.
\]
By \eqref{eq:projection-entropy-definition},
\[
  h_\pi(\sigma,\Pp)=\log N-R_1.
\]
Since $\Pp$ is ergodic and $\{S_a\}_{a\in\D}$ is a self-similar
IFS with common contraction ratio $\rho$, Feng and
Hu~\cite[Theorem~2.8]{FengHu2009} imply that $\mu$ is exact
dimensional and
\[
  \dimH\mu=\frac{h_\pi(\sigma,\Pp)}{-\log\rho}
  =\frac{\log N-R_1}{-\log\rho},
\]
which proves \eqref{eq:feng-hu-specialized}.
\end{proof}

\subsection{Proof of Theorem \ref{thm:rigidity} in Dimension One}

We now complete the proof of Theorem~\ref{thm:rigidity} in dimension one.

\begingroup
\setlength{\parskip}{0pt}
\setlength{\abovedisplayskip}{5pt plus 1pt minus 1pt}
\setlength{\belowdisplayskip}{5pt plus 1pt minus 1pt}
\begin{proof}[\bf Proof of Theorem~\ref{thm:rigidity} in dimension one]
Assertion (i) follows from Proposition~\ref{prop:dirac-fibers} and Corollary~\ref{cor:essential-injectivity}, which together characterize the essential injectivity of the coding projection and the almost-everywhere Dirac structure of the fiber measures.\par\noindent Corollary~\ref{cor:exact-entropy} establishes assertion (ii), namely that the fiber branching entropy vanishes for all $n\ge1$. Since $R_1=0$, Proposition~\ref{prop:feng-hu-interface} yields assertion (iii), namely
\[
  \dimH\mu=-\frac{\log N}{\log\rho}\le1.
\]
\end{proof}

In particular, a dimension drop
$\dimH\mu<\log N/(-\log\rho)$ is incompatible with the frame-measure
hypothesis.

\endgroup

\begin{remark}[Relation to Wang and Tang]\label{rem:wang-tang}
In dimension one, setting $d_a=b_a/\rho$ puts our maps in the form
$S_a(x)=\rho(x+d_a)$ considered by Wang and Tang~\cite{WangTang2026}.
In our notation, their conditional entropy $H(\omega_1\mid X)$ is
$R_1$, since $X=\pi(\omega)$. They proved that $R_1>0$ rules out a
frame measure. Hence frame-measure existence implies $R_1=0$. The
conditional-probability recursion and shift invariance then imply that
$\Pp_x$ is a Dirac mass for $\mu$-almost every $x$. Thus their result
gives the one-dimensional Dirac conclusion of
Theorem~\ref{thm:rigidity}. Our theorem also treats arbitrary ambient
dimensions.
\end{remark}

\section{Extension to Higher Dimensions}\label{sec:Rd}

We now prove the higher-dimensional version of Theorem~\ref{thm:rigidity}. The probabilistic and symbolic arguments developed after the collision estimate in Subsection~\ref{subsec:collision-estimate} do not rely on the order structure of $\R$.  It therefore remains only to replace the rightmost endpoint branch by a visible branch selected through an exposed translation, thereby verifying the common collision hypothesis \eqref{eq:mean-collision-lower}.

We first recall the canonical definition of an exposed point from convex geometry. 
A point $b$ in a compact convex set $C\subset\R^d$ is called
\emph{exposed} if there exists a vector $\ell\in\R^d$ satisfying
\[
  \ell\cdot b>\ell\cdot c, \qquad \forall c\in C\setminus\{b\}.
\]
All vertices of a polytope are exposed points.  Select an index $a_\ast$ such that the translation parameter
$b_{a_\ast}$ is a vertex of
$\operatorname{conv}\{b_0,\ldots,b_{N-1}\}$, and fix a separating direction $\ell\in\R^d$ with
\begin{equation}\label{eq:Rd-exposed}
  \ell\cdot b_{a_\ast}>\ell\cdot b_a,
  \qquad \forall a\ne a_\ast.
\end{equation}
Define the fixed point of the affine map \(S_{a_\ast}\) and its iterated cylinder set by
\[
  x_{a_\ast}:=\frac{b_{a_\ast}}{1-\rho},
  \qquad
  I_n:=S_{a_\ast}^{\circ n}(K).
\]

\vspace{0.2cm}

The next lemma gives the geometric separation needed for the collision estimate.

\begin{lemma}\label{lem:Rd-visible-branch}
The fixed point $x_{a_\ast}$ lies outside the image $S_a(K)$ for any
$a\ne a_\ast$.  Moreover, there exist $n_0\ge1$ and $C_0>0$ such that
\begin{equation}\label{eq:Rd-visible-mass}
  \mu(I_n)=C_0N^{-n},
  \qquad \forall n\ge n_0.
\end{equation}
\end{lemma}

\begin{proof}
Define the maximal linear functional value over the attractor \(K\) by
\[
  M:=\max_{x\in K} \ell\cdot x.
\]
By the self-similar invariance of \(K\) and the form of each affine branch \(S_a(x)=\rho x + b_a\), we rewrite the maximum via 
 \eqref{eq:Rd-exposed} as
\begin{align*}
  M
  &=\max_{a\in\D}\max_{x\in K}
       \bigl(\rho \ell\cdot x+\ell\cdot b_a\bigr)\\
  &=\rho M+\ell\cdot b_{a_\ast}.
\end{align*}
This gives that
\[
  M=\frac{\ell\cdot b_{a_\ast}}{1-\rho}
   =\ell\cdot x_{a_\ast}.
\]
Now take any $a\ne a_\ast$. For all \(y\in S_a(K)\), write \(y=S_a(x)=\rho x+b_a\) with \(x\in K\). We have
\[
  \max_{y\in S_a(K)}\ell\cdot y
  =\rho M+\ell\cdot b_a
  <\rho M+\ell\cdot b_{a_\ast}=M=\ell\cdot x_{a_\ast}.
\]
Hence $x_{a_\ast}\notin S_a(K)$.

We now construct the union of all non-\(a_\ast\) branch images:
\[
  F:=\bigcup_{a\ne a_\ast}S_a(K).
\]
The set $F$ and $\{x_{a_\ast}\}$ are disjoint and compact, so
$\delta:=\operatorname{dist}(\{x_{a_\ast}\}, F)>0$.  Since
\[
  S_{a_\ast}^{\circ n}(x)
  =x_{a_\ast}+\rho^n(x-x_{a_\ast}),
\]
the sets $I_n$ converge to $\{x_{a_\ast}\}$ in Hausdorff metric as \(n\to\infty\).
Choose $n_0$ sufficiently large so that $I_n\cap F=\varnothing$ for all $n\ge n_0$.

For any \(n\ge n_0\), \(I_n\subset S_{a_\ast}(K)\) and \(I_n\) intersects no other branch image. By the self-similar invariance relation for \(\mu\), the invariance identity gives the following recurrence
\[
  \mu(I_n)
  =\frac1N\mu(S_{a_\ast}^{-1}I_n)
  =\frac1N\mu(I_{n-1}).
\]
As $I_n\subset S_{a_\ast}(K)$, evaluating
\eqref{eq:intro-invariance} on $I_n$ leaves only the branch $a_\ast$.
Iterating this geometric recurrence yields \(\mu(I_n) = N^{-n+n_0}\mu(I_{n_0})\). Set \(C_0 := N^{n_0}\mu(I_{n_0})>0\), which proves equation \eqref{eq:Rd-visible-mass}. Positivity of \(C_0\) follows from the inclusion
$[a_\ast^{n_0}]\subset\pi^{-1}(I_{n_0})$, which guarantees
$\mu(I_{n_0})\ge N^{-n_0}>0$.
\end{proof}


\begin{proof}[\bf Proof of Theorem~\ref{thm:rigidity} in higher dimensions]
Recall that
\[
  c_n(x):=\sum_{\abs u=n} (p_u(x))^2.
\]
By Lemma~\ref{lem:posterior}, the measure
$\mu_{a_\ast^n}= p_{a_\ast^n}\mu=N^{-n}(S_{a_\ast}^{\circ n})_\#\mu$ is supported on
$I_n$ and has total mass $N^{-n}$. For $n\ge n_0$, Cauchy--Schwarz and
Lemma~\ref{lem:Rd-visible-branch} give
\[
  N^{-2n}
  =\left(\int_{I_n}p_{a_\ast^n}\dd\mu\right)^2
  \le\mu(I_n)\norm{p_{a_\ast^n}}_2^2
  =C_0N^{-n}\norm{p_{a_\ast^n}}_2^2.
\]
Thus
\[
  \norm{p_{a_\ast^n}}_2^2\ge C_0^{-1}N^{-n}.
\]
By Proposition~\ref{prop:comparability}, every word $u$ of length $n$
satisfies
\[
  \norm{p_u}_2^2\ge\frac{A}{BC_0}N^{-n}.
\]
Summing over all $N^n$ words with length $n$ yields
\begin{equation}\label{eq:Rd-collision-lower}
  \int_Kc_n(x)\dd\mu(x)
  =\sum_{\abs u=n}\norm{p_u}_2^2
  \ge\frac{A}{BC_0},
  \qquad n\ge n_0.
\end{equation}
Hence \eqref{eq:mean-collision-lower} holds with
$\kappa=A/(BC_0)$.

We now invoke the symbolic arguments of Sections~\ref{sec:atomicity} and~\ref{sec:entropy}. Their
proofs use only the disintegration over $\pi$, the refinement of cylinder
partitions, the Bernoulli shift, the injectivity of each branch, the coding
identity
\[
  \pi(a\tau)=S_a(\pi(\tau)),
\]
and the uniform collision bound.  All of these remain valid in $\R^d$.
For completeness, the core one-step fiber identity from
Proposition~\ref{prop:bayes} is
\[
  \Pp_x([a]\cap\sigma^{-1}(G))
  =\1_{S_a(K)}(x)r_a(x)
     \Pp_{S_a^{-1}(x)}(G),
\]
for $\mu$-almost every $x\in S_a(K)$, with the null exceptional set
independent of the Borel set $G\in\B(\Omega)$.  Therefore the arguments following
\eqref{eq:mean-collision-lower} in Sections~\ref{sec:atomicity} and~\ref{sec:entropy} force
the conditional coding measures $\Pp_x$ to be Dirac masses for
$\mu$-almost every $x\in K$.
The proof of Corollary~\ref{cor:essential-injectivity} is likewise
unchanged, and Corollary~\ref{cor:exact-entropy} gives
\[
  R_n=0,
  \qquad \forall n\ge1.
\]
Proposition~\ref{prop:feng-hu-interface} then yields exact dimensionality and
\[
  \dimH\mu=-\frac{\log N}{\log\rho}\le d.
\]
Consequently, $N\rho^d>1$ is impossible. This completes the proof of
Theorem~\ref{thm:rigidity} in arbitrary dimension.
\end{proof}

We conclude this section by recording the precise separation property provided by essential injectivity.

\begin{remark}\label{rem:essential-vs-osc}
Essential injectivity constitutes a measure-theoretic separation property. Analogous to the one-dimensional reasoning laid out after Corollary~\ref{cor:essential-injectivity}, this property yields mutual singularity for cylinder measures sharing the same word length. Importantly, it does not assume the stronger open set separation condition used in earlier results. Topological overlaps between branch images may still occur even under vanishing \(\mu\)-density. This inherent limitation stems from our \(L^2(\mu)\)-based analysis, which is insensitive to modifications carried out on null sets.
\end{remark}

\section{Concluding Remarks}\label{sec:further-questions}

Two questions remain. One concerns the parameters for which an $N$-Bernoulli convolution may admit a frame measure; the other is whether the equal-weight assumption can be relaxed.

\vspace{0.2cm}

For the uniform \(N\)-Bernoulli convolution
\(\mu_{\rho,N}\) defined in \eqref{eq:intro-bernoulli}, the critical (\(\rho =1/N\))
and supercritical regimes (\(\rho>1/N\)) can be described explicitly. If
\(\rho>1/N\), Corollary~\ref{cor:supercritical} rules out frame measures.
If \(\rho=1/N\), then
\[
  S_a([0,1])
  =\left[\frac{a}{N},\frac{a+1}{N}\right],
  \qquad a\in\D.
\]
These intervals cover \([0,1]\) with disjoint interiors. Hence the
normalized Lebesgue measure \(\1_{[0,1]}\dd x\) satisfies the
self-similar identity \eqref{eq:intro-invariance}, and uniqueness yields
\eqref{eq:intro-critical-measure}. Additionally,
\(\{e^{2\pi ikx}:k\in\Z\}\) is an orthonormal basis.  Consequently, for
 \(N^{-1}\le\rho<1\), the measure \(\mu_{\rho,N}\) admits a
frame measure if and only if \(\rho=N^{-1}\).

\vspace{0.2cm}

Within the separated regime $0<\rho<1/N$, several arithmetic
nonexistence results are known. Pont, Liehr and Taylor
~\cite{PontLiehrTaylor2026} proved that Cantor measures with odd
integer base do not admit Fourier frames. Fu, Song and Wang
~\cite{FuSongWang2026} obtained a Walsh--quotient obstruction for
certain reciprocal-power Bernoulli convolutions, and Fu and Wang
~\cite{FuWang2026} proved an integer-power obstruction for arbitrary
$N$. Wang and Tang~\cite{WangTang2026} extended this type of obstruction
to certain reduced rational powers: for consecutive digits, if
$\rho^{-m}=a/b$ in lowest terms, $b\mid N$, and $N/b\nmid a$, then
no frame measure exists. A full classification for the remaining
separated parameters is still open.

\vspace{0.2cm}

Our results rely on the equal-weight hypothesis. Given a probability weight vector  \(\mathbf p=(p_0,\ldots,p_{N-1})\) with \(p_a>0\) and \(\sum_{a=0}^{N-1}p_a=1\), define the weighted Bernoulli convolution
\[
  \mu_{\mathbf p}=\sum_{a=0}^{N-1}p_a(S_a)_\#\mu_{\mathbf p}.
\]
It is natural to ask whether, without a separation assumption, the existence of a frame measure forces \(p_0=\cdots=p_{N-1}\). More generally, does frame-measure existence force the essential injectivity of the coding map and the dimension identity
\[
  \dim_H\mu_{\mathbf p}=\frac{-\sum_{a=0}^{N-1}p_a\log p_a}{-\log\rho}?
\]
Our prior uniform comparison in Proposition \ref{prop:comparability} across identical word lengths fails to generalize directly.
For a finite word \(u=(u_1,\ldots,u_n)\), the associated cylinder measure now becomes
\[
  \mu_u=(p_{u_1}\cdots p_{u_n})(S_u)_\#\mu_{\mathbf p}.
\]
Its Fourier modulus depends sensitively on the exact symbol sequence of \(u\), rather than solely the word length. The main difficulty in extending the argument to nonuniform weights is to obtain uniform control of these cylinder contributions.

\begingroup
\fontsize{9.5}{10.5}\selectfont
\setlength{\parskip}{0pt}

\medskip
\noindent\textsc{Wei-Jie Wang}\\[-1pt]
Hubei Key Laboratory of Mathematical Sciences, College of Mathematics and
Statistics, Central China Normal University, Wuhan, Hubei 430079, China\\[-1pt]
Email addresses: 
\texttt{wwjmath@163.com} (W.-J. Wang)
\endgroup

\end{document}